\documentclass{amsart}
\usepackage{amssymb}
\usepackage{amsmath}
\usepackage{amscd}
\usepackage{graphics}
\usepackage{latexsym}

\def\C{\mathbb{C}}

\def\Z{\mathbb{Z}}

\newtheorem*{main}{Main Theorem}
\newtheorem{theorem}{Theorem}[section]
\newtheorem{lemma}[theorem]{Lemma}
\newtheorem{proposition}[theorem]{Proposition}

\newtheorem{example}[theorem]{Example}

\def\CP{\hbox{${\mathbb C}P^2$}}
\def\CPB{\hbox{$\overline{{\mathbb C}P^2}$}}
\def\HCP{\hbox{$H{\mathbb C} P^3$}}

\title{A New Construction of $6$-Manifolds}
\author{Ahmet Beyaz}
\address{Department of Mathematics UC Irvine CA 92697}
\email{abeyaz@math.uci.edu}
\address{Department of Mathematics Middle East Technical University, Ankara 06531 Turkey}
\email{beyaz@metu.edu.tr}
\thanks{}

\begin{document}

\begin{abstract}
This paper provides a topological method to construct all simply-connected, spin, smooth $6$-manifolds with torsion-free homology using simply-connected, smooth $4$-manifolds as building blocks. We explicitly determine the invariants that classify these $6$-manifolds from the intersection form and specific homology classes of the $4$-manifold building blocks. 
\end{abstract}

\maketitle
\setcounter{section}{-1}

\section{Introduction} \label{introduction}

The goal of this paper is to provide an explicit construction of all smooth, closed, simply-connected, and spin $6$-dimensional manifolds with torsion-free homology. A spin manifold is an oriented manifold such that the second Stiefel-Whitney class, $w_2$, of the tangent bundle is zero. We construct these $6$-dimensional manifolds using a plumbing construction on $2$-disk bundles over a carefully chosen collection of smooth, simply-connected $4$-dimensional manifolds. We then explicitly relate specific characteristic classes and invariants of these $4$-manifolds with the invariants that classify the $6$-manifold. In Section~\ref{4man}, we provide the existence results for the $4$-manifolds that we use in the constructions. 

Manifolds of dimension $6$ have been completely classified by C. T. C. Wall \cite{Wa-66} using standard algebraic topological invariants.

{\theorem [Wall] \label{6-man} Orientation-preserving diffeomorphism classes of simply-connected, smooth, spin, closed 6-manifolds $M$ with torsion-free (co)homology correspond bijectively to isomorphism classes of systems of invariants consisting of
\begin{enumerate}
\item a free Abelian group $H$ $(=H^2(M;\Z))$,
\item a symmetric trilinear map $\mu: H\times H\times H\rightarrow Z$ defined by $\mu(x,y,z)=x\cup y\cup z[M]$ satisfying
$\mu(x,x,y) \equiv \mu(x,y,y) \quad (\text{mod}\,2)$ for all $x,y \in H$,
\item a homomorphism $p_1: H\rightarrow Z$ such that
$p_1(x) \equiv 4\mu(x,x,x) \quad (\text{mod}\,24)$ for all $x \in H$.
\item a nonnegative integer $r = b_3/2$.
\end{enumerate}}

Wall proved Theorem~\ref{6-man} by using surgerical methods and homotopy information associated with these surgeries. In \cite{Wa-66}, Wall also proves that every $6$-manifold $M$ is diffeomorphic to $M_0\#_{b_3/2}S^3\times S^3$ with $M_0$ having the same invariants as $M$ but with $r=0$. The manifold $M_0$ is called the core of $M$. Therefore without loss of generality, our constructions deal with those $6$-manifolds with $r=0$.

An interesting family of simply-connected, spin $6$-manifolds are homotopy complex projective $3$-spaces. $6$-Manifolds which are homotopy equivalent to ${\mathbb C}P^3$ are called homotopy complex projective 3-spaces and each of them is generically denoted by \HCP. Their cohomology rings are isomorphic to that of ${\mathbb C}P^3$. The diffeomorphism types of \HCP's are distinguished by their first Pontrjagin classes and they are in one to one correspondence with the set of integers. In \cite{MY-66}, Montgomery and Yang give two characterizations of these manifolds. First they show that the diffeomorphism types of \HCP's are in one to one correspondence with free differentiable actions of $S^1$ on a homotopy $S^7$. It is also shown that smooth $S^1$-actions on a homotopy $S^7$ which have an $S^3$ as the fixed point set are in one to one correspondence with the isotopy classes of pairs $(S^6,K)$ where $K$ is an embedding of $S^3$ into $S^6$. These are exactly the smooth $3$-knots in $S^6$ (or so-called Haefliger knots, see \cite{Ha-62} and \cite{Ha-66}). The set of Haefliger knots is also in one to one correspondence with the integers. The second characterization of \HCP's is based on a surgery on these knots. If we remove a neighboorhood $B^3\times K^3$ of $K$ from $S^6$, the remaining part, the knot complement, is diffeomorphic to $S^2\times B^4$. Each \HCP~ can be formed by taking out a neighborhood $B^3\times S^3$ of the knot from $S^6$ and attaching a $S^2\times B^4$ to the knot complement. The attaching map $f$ from $S^2\times S^3$ to itself is different from the ordinary one, otherwise the space after the attaching would be $S^2\times S^4$. Montgomery and Yang do not give an explicit formula for $f$ except for $\mathbb{C}P^3$, but they point out that the action of $f$ on the third homotopy group of $S^2\times S^3$ must satisfy two conditions. Let $\rho_1:S^2\times S^3\rightarrow S^2$ and $\rho_2:S^2\times S^3\rightarrow S^3$ be the projection maps of $S^2\times S^3$ to $S^2$ and $S^3$, respectively. Let $\varphi:S^3\rightarrow S^2\times S^3$ be the inclusion map. Then the self-diffeomorphism $f$ of $S^2\times S^3$ must satisfy the following two conditions:

(1) $\rho_1\circ f\circ \varphi:S^3\rightarrow S^2$ represents a generator of $\pi_3(S^2)$. (i.e. $\rho_1\circ f\circ \varphi$ is the Hopf fibration)

(2) $\rho_2\circ f\circ \varphi:S^3\rightarrow S^3$ is of degree $-1$.

In Section~\ref{constructions}, we give a new way to construct \HCP's. This construction is a model for the later constructions of all $6$-manifolds. We take some $2$-disk bundles on the $4$-manifold and then close the boundary to get the \HCP.

When the second Betti number is greater than one, we need a surgery method which is called plumbing. Given two spheres $\Sigma_1$ and $\Sigma_2$, let $N_i$ be the total space of the direct sum of two $2$-disk bundles $\eta_{i1}$ and $\eta_{i2}$ over $\Sigma_i$. The plumbing on $\Sigma_1$ and $\Sigma_2$ is done by gluing $N_1$ and $N_2$ under a map that takes $\Sigma_1$ to $\Sigma_2$ and exchanges the factors of the direct sum pointwise. If $\Sigma_1$ and $\Sigma_2$ are embedded in $6$-manifolds $M_1$ and $M_2$, respectively, then plumbing of $M_1$ and $M_2$ on the spheres $\Sigma_1$ and $\Sigma_2$ is done by identifying the normal neighborhoods of the spheres in the respective manifold. A more detailed description is given in Section~\ref{constructions}. Now let us state the main theorem of the paper.

{\main \label{maintheorem} Let $V$ be a smooth, closed, simply-connected, spin $6$-manifold with torsion-free homology and $b_3(V)=0$. Suppose that $H^2(V;\Z)$ is isomorphic to the direct sum of $n$ copies of $\Z$, each of which is generated by $x_i$ ($1\leq i\leq n$). Also suppose that $p_1(V)x_i=24k_i+4\mu(x_i,x_i,x_i)$, where $p_1$ is the first Pontrjagin class of $V$ and $\mu$ is the symmetric trilinear form of $V$.
Then we can find a collection, $\{X_i\}$, of smooth, closed, simply-connected $4$-manifolds with odd intersection forms $Q_i$ and second cohomology classes $\alpha_{ij}\in H^2(X_i;\Z)$ ($1\leq i,j\leq n$) satisfying 

\begin{enumerate} \item the signature of $X_i$ is $8k_i+\mu(x_i,x_i,x_i)$,
\item $\alpha_{ii}$ are primitive, characteristic and $Q_i(\alpha_{ii},\alpha_{ii})=\mu(x_i,x_i,x_i)$,
\item if $i\neq j$ then $\alpha_{ij}$ has a smooth sphere representative in $X_i$,
\item $Q_i(\alpha_{ij},\alpha_{ik})=\mu(x_i,x_j,x_k)$,
\end{enumerate}

so that the manifold $M$ that is constructed by closing the boundaries of the plumbed $2$-disk bundles over these $4$-manifolds $X_i$ with Euler class $\alpha_{ii}$ is diffeomorphic to $V$. The plumbing of the respective bundles over $X_i$ is done on the sphere representatives of $\alpha_{ij}$ in $X_i$ and $\alpha_{ji}$ in $X_j$.}

Note that, $\alpha\in H^2(X;\Z)$ is characteristic if and only if $Q(\alpha,\beta)\equiv Q(\beta,\beta)$ for all $\beta\in H^2(X;\Z)$. If $H^2(X;\Z)$ has no torsion, $\alpha$ is primitive if and only if the subspace of $H^2(X;\Z)$ obtained by modding out the subspace generated by $\alpha$ has no torsion and its rank is $b_2(X)-1$.

The motivation behind these constructions is to study the symplectic structures on the $6$-manifolds and relate them to the smooth or symplectic topology of the $4$-manifolds. As a starting point, we observe that while all the exotic smooth structures on a given $4$-manifold vanishes when crossed with the $2$-sphere $S^2$, it still may be the case that the symplectic exoticness is retained. This is hinted at in the early work of Ruan \cite{Ru-94}. 

There is not much known about $6$-dimensional symplectic topology until now. There is no general method to distinguish the symplectic structures on a smooth $6$-manifold. Also it seems hard to decide whether a smooth $6$-manifold is symplectic or not. In particular, despite their  fairly simple topology, it is unknown which of the \HCP's possess a symplectic structure, other than ${\mathbb C} P^3$ itself. With the construction given in this paper, this problem is replaced with a problem of symplectic surgery. In order to use the results of this paper to study the symplectic structures on the $6$-manifolds considered in Subsections~\ref{b=2} and \ref{b=n}, first one must show that plumbing is a symplectic operation. Moreover, we must know which of the $4$-manifolds that we use as the building blocks can be chosen to be symplectic.

Before we start building the manifolds, we state a theorem about certain $5$-manifolds that appears as the boundaries of the $6$-manifolds in the intermediate steps of our constructions. 

{\theorem [Duan-Liang (\cite{DL-05})] \label{DL} Assume that $X$ is a closed, simply-connected, smooth $4$-manifold and $\alpha \in H^{2}(X;\Z)$ a primitive, characteristic class. If $X_\alpha$ is the total space of the $S^1$-bundle over $X$ with Euler class equal to $\alpha$, then $X_\alpha$ is diffeomorphic to $\#_{b_2(X)-1}S^2\times S^3$.}

\section{The Constructions} \label{constructions}

Our constructions use a certain collection of $4$-manifolds as the building blocks. We establish the existence of the appropriate $4$-manifolds in Section~\ref{4man}.

\subsection{Spin 6-Manifolds with $b_2=1$} \label{b=1}

Let $X$ be a closed, simply-connected, smooth $4$-manifold and $\alpha \in H^2(X)$ be a primitive, characteristic element. A $2$-disk bundle over $X$ is characterized by its Euler class $\alpha\in H^2(X;\Z)$. If $M_\alpha$ is the total space of the $2$-disk bundle $\zeta$ over $X$ with Euler class $\alpha$, then $\partial M_{\alpha}$ is a circle bundle over $X$ with its Euler class equal to $\alpha$. By Theorem~\ref{DL}, $\partial M_\alpha$ is diffeomorphic to $\#_{b_2(X)-1} S^2 \times S^3$. Denote the manifold that is constructed by attaching a $B^3 \times S^3$ to each component in the connected sum by $M$. More precisely, we get $M$ by capping off the boundary of $M_\alpha$ with $\natural_{b_2(X)-1}B^3 \times S^3$, the boundary sum of $b_2(X)-1$ copies of $B^3 \times S^3$. The invariants of the closed $6$-manifold $M$ are given in the following propositions.

{\proposition \label{sc} $M$ is simply-connected.}
\begin{proof} The homotopy equivalence between $M_\alpha$ and $X$ and the simply-connectedness of $X$ implies the simply-connectedness of $M_\alpha$. The boundary of $M_{\alpha}$ is diffeomorphic to $\#_{b_2(X)-1} S^2 \times S^3$ (Theorem~\ref{DL}), which is connected and simply-connected. Notice that $\natural_{b_2(X)-1}B^3 \times S^3$ is also simply-connected. Then by the Van-Kampen's theorem $M$ is simply-connected.
\end{proof}

{\proposition \label{betti1} For $M$, $b_0=b_2=b_4=b_6=1$ and $b_1=b_3=b_5=0$.}
\begin{proof} Since $M$ is connected, $H_0(M;\Z)=\Z$. The abelianization of the fundamental group of $M$ is $H_1(M;\Z)$. The fundamental group is trivial, so $H_1(M;\Z)$ is trivial. The manifold $M$ is closed and oriented, therefore $H_6(M;\Z)=\Z$ and Poincare duality can be applied. This gives $H^5(M;\Z)=0$. Since there are no torsion elements, $H_5(M;\Z)=0$ by the Universal Coefficient theorem. Similarly, $b_2(M)=b_4(M)$. Here is a part of the Mayer-Vietoris sequence that applies to the construction:

\begin{equation*}
0\rightarrow H_4(M_{\alpha};\Z)\oplus H_4(\natural_{b_2(X)-1}B^3 \times S^3;\Z) \rightarrow H_4(M;\Z)\rightarrow H_3(\#_{b_2(X)-1}S^2\times S^3;\Z)
\end{equation*}

\begin{equation} \label{MV1}
\rightarrow H_3(M_{\alpha};\Z)\oplus H_3(\natural_{b_2(X)-1}B^3 \times S^3;\Z) \rightarrow H_3(M;\Z) \rightarrow H_2(\#_{b_2(X)-1}S^2\times S^3;\Z)
\end{equation}

\begin{equation*}
\rightarrow H_2(M_{\alpha};\Z)\oplus H_2(\natural_{b_2(X)-1}B^3 \times S^3;\Z) \rightarrow H_2(M;\Z) \rightarrow 0
\end{equation*}

By the homotopy equivalence of $M_{\alpha}$ and $X$, $H_4(M_{\alpha};\Z)=\Z$, and by the homotopy equivalence of $\natural_{b_2(X)-1}B^3 \times S^3$ and a wedge of $S^3$'s, $H_4(\natural_{b_2(X)-1}B^3 \times S^3)=0$. Because of this fact, the map below which is induced by the inclusion map is an isomorphism.
$$H_3(\#_{b_2(X)-1}S^2\times S^3;\Z) \rightarrow H_3(M_{\alpha};\Z)\oplus H_3(\natural_{b_2(X)-1}B^3 \times S^3;\Z)$$
Therefore, exact sequence (\ref{MV1}) splits into three short exact sequences. First of them gives $b_4(M)=1$ and $b_2(M)=1$. The second part of the split sequence implies that $H_3(M;\Z)=0$.
\end{proof}

{\proposition \label{generator} $H^2(M;\Z)=\Z$ has a generator whose pullback in $X$ is $\alpha$.}
\begin{proof} The boundary of $M_{\alpha}$ is the circle bundle over $X$ with Euler class $\alpha$. A part of the Gysin sequence for this circle bundle over $X$ is given below.

\begin{equation} \label{Gysin}
0\overset{\pi^*}{\rightarrow }H^{3}(\#_{b_2(X)-1}S^2\times S^3;\Z)\rightarrow H^2(X;\Z)\overset{\cup \alpha }{\rightarrow }H^{4}(X;\Z) \overset{\pi^*}{\rightarrow} 0
\end{equation}

Kernel of $\cup\alpha$ is equal to the image of the second map which is injective. The image is the quotient of $H^2(X;\Z)$ by the subspace generated by a class such that intersection of this class with $\alpha$ is one, because $\cup\alpha$ is surjective. The elements in the second homology of $\#_{b_2(X)-1}S^2\times S^3$, given by the Poincare duals of the third cohomology are capped with $B^3$'s of $\natural_{b_2(X)-1}B^3 \times S^3$. The integer multiples of $\alpha$ are among the surviving cohomology elements. The injection of $X$ into $M$ induces the pullback of the generator $x$ to $X$ as $\alpha$.
\end{proof}

{\proposition \label{spinb=1} $M$ is spin.}
\begin{proof} It is enough to calculate the value of $w_2(M)$ on the generator $PD(x)$ of $H_2(M;\Z)=\Z$. Since inclusion of $H_2(M_\alpha)$ to $H_2(M)$ is surjective, $w_2(M)=j^*(w_2(M_{\alpha}))$, where $j$ is the inclusion map of $M_{\alpha}$ into $M$. Let $i$ be the inclusion map of $X$ into $M_{\alpha}$. Writing $i^*(TM_{\alpha})=TX \oplus \zeta$, the Whitney sum formula gives, $i^* w_2(M_{\alpha}) \equiv w_2(X)+\alpha$ (mod 2) in $H^2(X;\Z_2)$. In $H_4(M;\Z)$, $PD(x)=i_*PD(\alpha)$. Therefore, $w_2(M)PD(x)\equiv (i^* w_2(M_{\alpha}))PD(x) \equiv (w_2(X)+\alpha)PD(\alpha) \equiv 0$ mod $2$.
\end{proof}

{\proposition \label{intersectionb=1} $Q(\alpha,\alpha)=\mu(x,x,x)$ where $x$ is the generator of $H^2(M;\Z)$.}
\begin{proof} The Poincare dual of $x\in H^2(M;\Z)$ is $[i_*X]$, where $i$ is the inclusion map. In $H^2(X;\Z)$, $\alpha=i^*x$. $[X]=i^* PD(x)$ in $H_4(X;\Z)$. $\mu(x,x,x)=(x\cup x\cup x)[M]=(x\cup x)(x\cap [M])=(x\cup x)[PD(x)]=(x\cup x)[i_*X]=(i^*\alpha \cup i^*\alpha)[i_*X]=(\alpha\cup \alpha)[X]=Q(\alpha,\alpha)$.
\end{proof}

{\proposition \label{pontrjaginb=1} $p_1(M)=p_1(X)+\alpha\cup\alpha$, where $p_1$ is the first Pontrjagin class of $M$.}
\begin{proof}
To see this equality, let's consider the Whitney sum formula for Pontrjagin classes. The first Pontrjagin class of $M$, $p_1(M)$ is given as the second Chern class $c_2(TM\otimes\C)$ of the complexification of the tangent bundle of $M$. Since there is no contribution to the fourth cohomology of $M$ from $\natural_{b_2(X)-1}B^3 \times S^3$, it is safe to make the calculations in $M_{\alpha}$. The Whitney sum formula for the second Chern number of $TM_{\alpha}\otimes\C$ is $c_2(TM_{\alpha}\otimes\C)=c_2(TX\otimes\C) + c_1(TX\otimes\C)c_1(\nu X|_{M_\alpha}\otimes\C) + c_2(\nu X|_{M_\alpha}\otimes\C)$. The term which is the product of the first Chern classes has order two in the cohomology (\cite{MS-75} p.175). The cohomology of $X$ is torsion-free, thus the term consisting of the product of $c_1$'s is $0$. The second Chern class of the complexification of the $2$-disk bundle with the Euler class $\alpha$ (\cite{MS-75}). Consequently, $p_1(M)x=p_1(M)(PD(x))=p_1(X)(PD(x))+\alpha\cup\alpha(PD(x))=p_1(X)[X]+Q(\alpha,\alpha)$. By the Hirzebruch signature theorem for $4$-manifolds, this is equal to $3\sigma(X)+m$. $\sigma(X)=8k+m$ implies that $p_1(M)x=24k+4m$.
\end{proof}


Using the propositions given above, we prove the theorem below. This theorem is a special case of the main theorem where $b_2=1$.

{\theorem Let $V$ be any closed, simply-connected, spin $6$-manifold with torsion-free homology and $b_3(V)=0$. Assume that $H^2(V;\Z)$ is isomorphic to $\Z$ which is generated by $x$. Also assume that $p_1(V)x= 24k+4\mu(x,x,x)$, where $p_1$ is the first Pontrjagin class of $V$ and $\mu$ is the symmetric trilinear form of $V$. Then we can find a smooth, closed, simply-connected $4$-manifold $X$ with an odd intersection form $Q$ and a second cohomology class $\alpha\in H^2(X;\Z)$ satisfying

\begin{enumerate} \item the signature of $X$ is $8k+\mu(x,x,x)$,
\item $\alpha$ is primitive, characteristic and $Q(\alpha,\alpha)=\mu(x,x,x)$,
\end{enumerate}

so that the manifold $M$ that is constructed by the method described above is diffeomorphic to $V$.}

\begin{proof} The existence of an appropriate $4$-manifold is established below in by Lemma~\ref{8k+mExistence}. Hence, to prove the theorem, we must show that the constructed manifold $M$ has the same invariants as $V$. 

The manifold $M$ is simply connected by Proposition~\ref{sc}. The manifold $M$ is spin by Proposition~\ref{spinb=1}. As shown above in Proposition~\ref{betti1}, the Betti numbers are same as of $M$. Since $Q(\alpha,\alpha)=\mu(x,x,x)$, cohomology rings are isomorphic by Proposition~\ref{intersectionb=1}. The first Pontrjagin classes of $M$ and $V$ coincide by Proposition~\ref{pontrjaginb=1}. By the Wall's theorem (Theorem~\ref{6-man}) $M$ is diffeomorphic to $V$. 
\end{proof}

{\example \label{hcp3} \rm The family of homotopy complex projective $3$-spaces form an interesting collection of $6$-manifolds. Some information about these manifolds can be found in the introduction. \HCP's are distinguished by their first Pontrjagin classes and parametrized by $\Z$. Namely, for all $k\in\Z$, the corresponding $\HCP$ has $p_1(\HCP)$ evaluated on the second cohomology element is $24k+4$. Take $X$ such that $\sigma(X)=8k+1$ and $\alpha \in H^2(X;\Z)$ such that $\alpha\cup\alpha=1$. When we apply the construction, we end up with a $\HCP$ with its first Pontrjagin class evaluated on $x$ equal to $24k+4$. For example, if $X$ is $\CP$ and $\alpha$ is $h$ where $h$ generates $H^2(\CP;\Z)$, we get $\mathbb{C}P^3$. Here $b_2-1=0$, hence $M_{\alpha}$ is the Hopf fibration over $\CP$ and $\partial M_{\alpha}$ is $S^5$. Gluing the $6$-disk to the boundary gives $\mathbb{C}P^3$.}

{\example \label{quintic} \rm The smooth quintic hypersurface $Q$ in $\mathbb{C}P^4$, which is the zero set of a degree $5$ polynomial, has $b_3=204$. $Q_0$, the core of $Q$, has its second cohomology group $H^2(Q_0;\Z)$ isomorphic to $\Z$. Let's denote the generator of this group by $L$. $\mu(L,L,L)=5$ and $p_1(Q_0)L=-100$. We can construct $Q_0$ by the construction given above by taking a $4$-manifold with signature equal to $-35$ and $\alpha$ with self intersection $5$. For example we can choose $X$ as $\CP\#_{36}\CPB$ and $\alpha$ as $7h+3e_1+3e_2+\Sigma_3^{36}e_i$. This choice is not unique. In fact there are infinitely many choices. Another choice for $X$ is the degree $5$ hypersurface in $\mathbb{C}P^3$ which is diffeomorphic to $\#_9\CP\#_{44}\CPB$.}

\subsection{Spin $6$-Manifolds with $b_2=2$} \label{b=2}

In this subsection, we apply a construction similar to that in Subsection~\ref{b=1} in order to get the simply-connected, spin, torsion-free, smooth $6$-manifolds with $b_2=2$ and $b_3=0$. Let $X_i^4$, $(i=1,2)$ be closed, simply-connected, smooth $4$-manifolds, $Q_i$ be their respective intersection forms and $\alpha_{ij} \in H^2(X_i,\Z)$ be primitive elements where $j=1,2$ and $\alpha_{ii}$ are characteristic. Let $M_i$ be the total space of the $2$-disk bundle $\zeta_i$ over $X_i$ with Euler class $\alpha_{ii}$. Suppose that the Poincare dual of each cohomology class, $\alpha_{ij} (i\neq j)$, can be represented by an embedded sphere $\Sigma_{ij}$ in $X_i$. Our $6$-manifold $M$ is formed by gluing $M_1$ and $M_2$ over a neighborhood of the $\Sigma_{ij}$ in $M_i$  by a diffeomorphism explained below and then capping off the boundary. 

Let $\nu_{ij}$ be the normal $2$-disk bundle of the sphere $\Sigma_{ij}$ in the $4$-manifold $X_i$, and let $\zeta_{ij}$, with total space $N_{ij}$ be the restriction of $\zeta_{i}$ to the total space of $\nu_{ij}$. The attaching of $M_1$ and $M_2$ is along the total spaces $N_{12}$ and  $N_{21}$. Since $N_{ij}$ is tubular neighborhood of $\Sigma_{ij}$ in $M_i$, it is also the total space of the $B^4$-bundle $\eta_{ij}$ over $\Sigma_{ij}$.

To glue $M_1$ and $M_2$ along $N_{12}$ and $N_{21}$, we need them to be diffeomorphic. Since $\pi_2(SO(3))=\Z_2$, there are only two topologically distinct $4$-disk bundles over $S^2$, determined by their second Stiefel-Whitney class. Thus, for the neighborhoods to be diffeomorphic, we need $w_2(\eta_{12})=w_2(\eta_{21})$. Since the Euler class of $\zeta_i$ is $\alpha_{ii}$, we have that $w_2(\eta_{ij})=Q_i(\alpha_{ij},\alpha_{ii})+Q_i(\alpha_{ij},\alpha_{ij})\mod{2}$. However, since $\alpha_{ii}$ are characteristic, we have that $w_2(\eta_{ij})=0$. This implies that the bundles $\eta_{ij}$ are trivial and $N_{ij}$ are diffeomorphic to $S^2\times B^4$.

To perform our plumbing construction, we need a stronger condition. The normal bundle $\eta_{ij}$ of the sphere $\Sigma_{ij}$ in $M_i$ can be decomposed as the direct sum of two $2$-disk bundles (\cite{Be-96}). The first is the normal bundle $\nu_{ij}$ of the sphere $\Sigma_{ij}$ in the $4$-manifold $X_i$ and the second is the restriction of the $2$-disk bundle $\zeta_i$ over $X_i$ restricted to $\Sigma_{ij}$. In other words, $\eta_{ij}=\,\nu_{ij}\,\oplus\,\zeta_{i}|_{\Sigma_{ij}}$.  

In order to mimic the plumbing of $2$-disk bundles over surfaces, we glue $N_{ij}$ by a diffeomorphism switching the first and the second bundles of the decomposition of $\eta_{ij}$ over $\Sigma_{ij}$ with the total space $N_{ij}$. That is, we identify normal bundle $\nu_{12}$ of $\Sigma_{12}$ in $X_1$ to the restriction of the $2$-disk bundle $\zeta_2$ over $X_2$ to $\Sigma_{21}$ and vice versa. All the identifications are made in an orientation preserving way. In order to do this, however, we need the $2$-disk bundles that are identified to be the homotopic. This then places two conditions on the cohomology classes $\alpha_{ij}$. These conditions are $Q_1(\alpha_{11},\alpha_{12})=Q_2(\alpha_{21},\alpha_{21})$ and $Q_1(\alpha_{12},\alpha_{12})=Q_2(\alpha_{21},\alpha_{22})$. Therefore, identifying $N_{12}$ in $M_1$ and $N_{21}$ in $M_2$ by pointwise identification of the fibers of $\nu_{12}$ with $\zeta_2|_{\Sigma_{21}}$ and $\nu_{21}$ with $\zeta_1|_{\Sigma_{12}}$ is a well-defined operation. Let us denote the product manifold by $M_{12}$.

{\proposition \label{sc2a} $M_{12}$ is simply-connected.}
\begin{proof} The homotopy equivalence between $M_i$ and $X_i$ and the simply-connectedness of $X_i$ implies the simply-connectedness of $M_i$. Note that, $S^2\times B^4$ is also simply-connected. Since the gluing space for the plumbing operation is path-connected, by the Van-Kampen's theorem, $M_{12}$ is simply-connected.
\end{proof}

{\proposition \label{homology2a} For $M_{12}$, $b_0=1$, $b_2=b_2(X_1)+b_2(X_2)-1$, $b_4=2$ and $b_1=b_3=b_5=b_6=0$.}
\begin{proof} We know that $M_{12}$ is connected, so $H_0(M_{12};\Z)=\Z$. The abelianization of the fundamental group is $H_1(M_{12};\Z)$. The fundamental group is trivial, so $H_1(M_{12};\Z)$ is trivial. Writing the Mayer-Vietoris sequence of this step enables us to determine the Betti numbers.

\begin{equation*}
0\rightarrow H_6(M_1;\Z)\oplus H_6(M_2;\Z) \rightarrow H_6(M_{12};\Z)\rightarrow H_5(S^2 \times B^4;\Z)
\end{equation*}

\begin{equation*}
\rightarrow H_5(M_1;\Z)\oplus H_5(M_2;\Z) \rightarrow H_5(M_{12};\Z)\rightarrow H_4(S^2 \times B^4;\Z)
\end{equation*}

\begin{equation} \label{MV2}
\rightarrow H_4(M_1;\Z)\oplus H_4(M_2;\Z) \rightarrow H_4(M_{12};\Z)\rightarrow H_3(S^2 \times B^4;\Z)
\end{equation}

\begin{equation*}
\rightarrow H_3(M_1;\Z)\oplus H_3(M_2;\Z) \rightarrow H_3(M_{12};\Z) \rightarrow H_2(S^2 \times B^4;\Z)
\end{equation*}

\begin{equation*}
\rightarrow H_2(M_1;\Z)\oplus H_2(M_2;\Z) \rightarrow H_2(M_{12};\Z) \rightarrow 0
\end{equation*}

When we place the known values to sequence (\ref{MV2}), it is straightforward to get $b_2(M_{12})=b_2(X_1)+b_2(X_2)-1$, $b_4(M_{12})=2$, $b_0(M_{12})=1$ and $b_1(M_{12})=b_3(M_{12})=b_5(M_{12})=b_6(M_{12})=0$.
\end{proof}

By Theorem~\ref{DL}, $\partial M_i$ is diffeomorphic to $\#_{b_2(X_i)-1} S^2 \times S^3$. We now claim that the boundary of $M_{12}$ is again diffeomorphic to a connected sum of a number of $S^2\times S^3$'s. More precisely, $\partial M_{12}$ is diffeomorphic to $\#_{b_2(X_1)+b_2(X_2)-3}S^2\times S^3$. To see this, we write the relative homology sequence for the pair $(M_{12},\partial M_{12})$.

\begin{equation*}
0\rightarrow H_5(M_{12};\Z)\rightarrow H_5(M_{12},\partial M_{12};\Z)\rightarrow H_4(\partial M_{12};\Z)
\end{equation*}

\begin{equation*}
\rightarrow H_4(M_{12};\Z)\rightarrow H_4(M_{12},\partial M_{12};\Z)\rightarrow H_3(\partial M_{12};\Z)
\end{equation*}

\begin{equation} \label{RMV1}
\rightarrow H_3(M_{12};\Z)\rightarrow H_3(M_{12},\partial M_{12};\Z)\rightarrow H_2(\partial M_{12};\Z)
\end{equation}

\begin{equation*}
\rightarrow H_2(M_{12};\Z)\rightarrow H_2(M_{12},\partial M_{12};\Z)\rightarrow H_1(\partial M_{12};\Z)
\end{equation*}

\begin{equation*}
\rightarrow H_1(M_{12};\Z)\rightarrow H_1(M_{12},\partial M_{12};\Z)\rightarrow \widetilde{H}_0(\partial M_{12};\Z)=0
\end{equation*}

By Poincare-Lefschetz duality and the Universal Coefficient theorem, we see that $H_2(M_{12},\partial M_{12})$ is isomorphic to $\Z\oplus\Z$ and $H_3(M_{12},\partial M_{12})=H_3(M_{12})$. The fourth homology of $M_{12}$ is consisting of the manifolds $X_1$ and $X_2$, so the map $H_4(\partial M_{12};\Z)\rightarrow H_4(M_{12};\Z)$ is the zero map. Now, using Poincare duality for $\partial M_{12}$, it is clear that the primitive classes that we are using for the plumbing are reflected as essential second homology classes in the boundary. By using sequence (\ref{RMV1}), we find the Betti numbers of $\partial M_{12}$ as $b_0=b_5=1$, $b_1=b_4=0$ and $b_2=b_3=b_2(X_1)+b_2(X_2)-3$. By the classification of simply-connected, spin, smooth $5$-manifolds (\cite{Sm-62}), the boundary is diffeomorphic to $\#_{b_2(X_1)+b_2(X_2)-3}S^2\times S^3$.

Once we know the boundary is diffeomorphic to $\#_{b_2(X_1)+b_2(X_2)-3}S^2\times S^3$, it is easy to obtain the desired manifold $M$ by capping off the boundary of $M_{12}$ with $\natural_{b_2(X_1)+b_2(X_2)-3}B^3\times S^3$. The invariants of $M$ are calculated in the following propositions.

{\proposition \label{sc2} $M$ is simply-connected.}
\begin{proof} The manifold $\natural_{b_2(X_1)+b_2(X_2)-3}B^3 \times S^3$ is simply-connected. By Proposition~\ref{sc2a}, $M_{12}$ is also simply-connected. Since the gluing space for the operation, $\#_{b_2(X_1)+b_2(X_2)-3}S^2\times S^3$, is path-connected, by the Van-Kampen's theorem, $M$ is simply-connected. 
\end{proof}

{\proposition \label{homology2} For $M$, $b_0=b_6=1$, $b_4=b_2=2$ and $b_1=b_3=b_5=0$.}
\begin{proof} $M$ is connected, so $H_0(M;\Z)=\Z$. The abelianization of the fundamental group is $H_1(M;\Z)$. The fundamental group is trivial, so $H_1(M;\Z)$ is trivial. The manifold $M$ is closed oriented, therefore $H_6(M;\Z)=\Z$ and Poincare duality can be applied. This gives $H^5(M;\Z)=0$. Since there is no torsion elements, by the Universal Coefficient theorem, $H_5(M;\Z)=0$. Similarly, $b_2(M)=b_4(M)$. 

The calculation of the remaining Betti numbers of $M$ is done by writing the Mayer-Vietoris sequence of the steps in the construction. There are two steps. 

The first step of the construction is attaching $M_1$ to $M_2$ by the map explained before Proposition~\ref{sc2a}, the one that plumbs $\alpha_{12}$ to $\alpha_{21}$. In Proposition~\ref{homology2a}, the Betti numbers of $M_{12}$ are given as $b_2(M_{12})=b_2(X_1)+b_2(X_2)-1$, $b_4(M_{12})=2$, $b_0(M_{12})=1$ and $b_1(M_{12})=b_3(M_{12})=b_5(M_{12})=b_6(M_{12})=0$.

The second step is closing the boundary of $M_{12}$. We calculate the Betti numbers as $b_0(M)=b_6(M)=1$, $b_1(M)=b_3(M)=b_5(M)=0$, and $b_2(M)=b_4(M)=2$ by using the following Mayer-Vietoris sequence.

\begin{equation*}
0\rightarrow H_4(M_{12};\Z)\oplus H_4(\natural_{b_2(X_1)+b_2(X_2)-3}B^3 \times S^3;\Z)\rightarrow H_4(M;\Z)
\end{equation*}

\begin{equation} \label{MV3}
\rightarrow H_3(\#_{b_2(X_1)+b_2(X_2)-3}S^2\times S^3;\Z)\rightarrow H_3(M_{12};\Z)\oplus H_3(\natural_{b_2(X_1)+b_2(X_2)-3}B^3 \times S^3;\Z)
\end{equation}

\begin{equation*}
\rightarrow H_3(M;\Z) \rightarrow H_2(\#_{b_2(X_1)+b_2(X_2)-3}S^2\times S^3;\Z)
\end{equation*}

\begin{equation*}
\rightarrow H_2(M_{12};\Z)\oplus H_2(\natural_{b_2(X_1)+b_2(X_2)-3}B^3 \times S^3;\Z) \rightarrow H_2(M;\Z) \rightarrow 0
\end{equation*}
\end{proof}

Next proposition shows that $M$ is spin.

{\proposition \label{spinb=2} $M$ is spin.}
\begin{proof} The proof is similar to the proof of Proposition~\ref{spinb=1}, but in this case there are two homology classes on which $w_2(M)$ must be evaluated. Both of these cohomology classes are induced by the inclusion of $X_i$ into $M$. Pulling back everything into respective $X_i$, it turns out that $w_2(M)x_i\equiv w_2(X_i)x_i+w_2(\zeta_i)x_i$. Remember that $\zeta_i$ is the $2$-disk bundle over the $4$-manifold $X_i$. In this sum, $w_2(X_j)x_i$ and $w_2(\zeta_j)x_i$ $(i\neq j)$ are counted due to the fact that $X_i$ overlaps $X_j$ only on $\alpha_{ij}$ which is already included in $x_i$. Hence $w_2(X_i)x_i$ and $w_2(\zeta_2)x_i$ are already included in $w_2(X_1)x_i$ and $w_2(\zeta_1)x_i$, respectively, and they are omitted. Therefore, $w_2(M)x_i\equiv w_2(X_i)x_i+w_2(\zeta_i)x_i\equiv w_2(X_i)x_i+Q_i(\alpha_{ii},\alpha_{ii}) \equiv 0$.
\end{proof}

As calculated in Proposition~\ref{homology2}, $H^2(M;\Z)$ is isomorphic to $\Z\oplus\Z$. This group is generated by the pushforwards of the cohomology elements $\alpha_{11}$ and $\alpha_{22}$ by the inclusion maps. We denote these two generators $x_1$ and $x_2$, respectively. The intersection form is determined by the values of $\mu(x_1,x_1,x_1)$, $\mu(x_1,x_1,x_2)$, $\mu(x_1,x_2,x_2)$ and $\mu(x_2,x_2,x_2)$. The embedding map of $X_i$ into $M$ is given by the composition of the inclusion maps of $f_i:X_i\rightarrow M_i$ and $g_i:M_{12}\rightarrow M$. Using the exact sequence (\ref{MV2}) above, we see that the inclusions $H^4(M_i)\rightarrow H^4(M_{12})$ are injective. The strategy for calculating the intersection numbers of $M$ is to focus on the intersection of these representatives in $M$. By transversality, we see that $PD(x_1)$ and $PD(x_2)$ in $H_4(M;\Z)$ intersect each other (and themselves) in a $2$-dimensional subspace.

Note that, the representatives of the Poincare duals to the second cohomology classes $x_1$ and $x_2$ are the embedded copies of the $4$-manifolds $X_1$ and $X_2$, respectively. This can be seen from the following argument. The Poincare dual of $x_i$, $PD(x_1)$, is in $H_4(M;,\Z)$. From the homology calculations above, $H_4(M;,\Z)$ is generated by $X_1$ and $X_2$. This means that $PD(x_1)$ is a linear combination of $X_1$ and $X_2$, say $sX_1+tX_2$. The normal bundle of $X_1$ in $M$ has the Euler class $\alpha_{11}$, which lies completely in $X_1$. The second cohomology class $x_1$ is the image of $\alpha_{11}$ in $M$, therefore the evaluation of $x_1$ on $X_2$ is restricted to $X_1\pitchfork X_2$ (the transversal intersection of $X_1$ and $X_2$ in $M$). This evaluation occurs within $X_1$, hence $t=0$. Since $PD(x_1)$ is primitive (due to the fact that $x_1$ is primitive in the second cohomology), $s=1$ and $PD(x_1)=X_1$. Similarly, $PD(x_2)=X_2$.

Our construction clearly gives that $PD(x_1) \pitchfork PD(x_2)=X_1\pitchfork X_2$ is the representing (smooth) surface of $\alpha_{12}$ in $X_1$ and the representing (smooth) surface of $\alpha_{21}$ in $X_2$. Moreover, $X_i$ has the normal bundle with the Euler class $\alpha_{ii}$ in $M$. As one of the consequences of the Thom isomorphism theorem (\cite{BT-82}, \cite{Br-93} p.382), $X_i$ intersects itself on the surface representing $\alpha_{ii}$. From these considerations, we may conclude the following proposition.

{\proposition Let $M$ be the $6$-manifold constructed above with $b_2(M)=2$. Suppose that $H^2(M;\Z)$ is generated by $x_1$ and $x_2$. Then $\mu(x_1,x_1,x_1)=Q_1(\alpha_{11},\alpha_{11})$, $\mu(x_1,x_1,x_2)=Q_1(\alpha_{11},\alpha_{12})$, $\mu(x_1,x_2,x_2)=Q_2(\alpha_{21},\alpha_{22})$ and $\mu(x_2,x_2,x_2)=Q_2(\alpha_{22},\alpha_{22})$.}
\begin{proof} 
\begin{enumerate} \item $\mu(x_i,x_i,x_i)=PD(x_i) \pitchfork PD(x_i) \pitchfork PD(x_i)=X_i\pitchfork X_i\pitchfork X_i$. The last $X_i$ intersects each of the other two on $\alpha_{ii}$ therefore the intersection of the three is equal to $Q_1(\alpha_{ii},\alpha_{ii})$.
\item $\mu(x_i,x_i,x_j)=PD(x_i) \pitchfork PD(x_i) \pitchfork PD(x_j)=X_i\pitchfork X_i\pitchfork X_j$. $(i\neq j$.) The first $X_i$ intersects the second $X_i$ on $\alpha_{ii}$ and intersects $X_j$ on $\alpha_{ij}$. As a result, these three intersect each other on the projection of $\alpha_{ij}$ on $\alpha_{ii}$ in $X_i$, which is nothing but $Q_i(\alpha_{ii},\alpha_{ij})$.
\end{enumerate}\end{proof}

Since $\mu$ is a symmetric trilinear form, this proposition determines it completely (\cite{OV-95}). The last piece of information we need to know about $M$ is the linear form on the second cohomology (i.e the first Pontrjagin class evaluated on the second cohomology elements) which we calculate now. 

{\proposition \label{pontrjagin2} $p_1(M)x_i=3\sigma(X_i)+Q_i(\alpha_{ii},\alpha_{ii})$.}
\begin{proof} We may write $p_1(M)x_i=p_1(M)PD(x_i)=p_1(M)X_i$. The tangent bundle of $M$ is restricted to $X_i$ as the direct sum of the tangent bundle $TX$ of $X$ and the normal neighborhood $\nu X_i|_M$ of $X_i$ in $M$. The bundle $\nu X_i|_M$ is the $B^2$-bundle $\zeta_i$ over $X_i$ with Euler class $\alpha_{ii}$. Consequently, as we have seen in Subsection~\ref{b=1}, $p_1(M)X_i=p_1(X_i)[X_i]+\alpha_{ii}\cup\alpha_{ii}[X_i]=3\sigma(X_i)+Q_i(\alpha_{ii},\alpha_{ii})$.
\end{proof}

\subsection{Spin $6$-Manifolds with $b_2=n$} \label{b=n}

Using the methods of the previous two subsections, it is possible to construct all simply-connected, spin, torsion-free, smooth $6$-manifolds with arbitrary second betti number and arbitrary trilinear form. (Note: There are some admissibility conditions for the trilinear forms to be realized by $6$-manifolds. Here, we are not going to study them. See \cite{OV-95} and \cite{Sc-97}.) Once we choose our $4$-manifolds with the suitable second homology groups and the intersection forms, whose existence is guaranteed by Lemma~\ref{b2=nExistence}, the construction is quite similar to the $b_2=2$ case. However, not all the manifolds in the intermediate steps are simply-connected. There are some $1$-dimensional homology elements that are created in the process, so we focus on the formation of these elements.

To construct a $6$-manifold $M$ with $b_2(M)=n$, first we take $n$ $4$-manifolds $X_i$ and $n$ primitive characteristic second cohomology classes $\alpha_{ii}\in H^2(X_i;\Z)$. We also take primitive second cohomology classes $\alpha_{ij}\in H^2(X_i;\Z)$ ($1\leq i,j \leq n,\, i\neq j$) each of which is represented by an embedded sphere. Let $M_i$ be the total space of the $2$-disk bundle $\zeta_i$ over $X_i$ with the Euler class $\alpha_{ii}$. Attach $M_i$ to $M_j$ by plumbing the image of the sphere representative of $\alpha_{ij}$ in $M_i$ under the inclusion map of $X_i$ on the image of the sphere representative of $\alpha_{ji}$ in $M_j$ under the inclusion map of $X_j$. Since the sphere representatives in any $X_i$ may not be disjoint (in fact, we need them to have nonempty intersection in general), all the plumbing must be done one by one.

We form manifolds $M_{ij}$ inductively by doing the plumbing one by one. Here $i$ and $j$ must be as if $ij$ is an index for the terms of a matrix above the diagonal. Let $M_{12}$ be the manifold obtained by attaching $M_1$ to $M_2$ by gluing $\alpha_{12}$ on $\alpha_{21}$ as explained in Subsection~\ref{b=2}. Similarly, $M_{13}$ is the manifold obtained by attaching $M_3$ to $M_{12}$ on $\alpha_{13}$ and $\alpha_{31}$. The manifold $M_{23}$ is the manifold obtained by attaching $M_3$ to $M_{1n}$ on $\alpha_{23}$ and $\alpha_{32}$. Then $M_{ij}$ is the manifold obtained by plumbing $M_j$ to the manifold formed in the last step ($M_{i(j-1)}$ or $M_{(i-1)j}$) on $\alpha_{ij}$ and $\alpha_{ji}$. $\frac{n(n-1)}{2}$ manifolds are formed this way. The last one is $M_{(n-1)n}$.

It is similar to the last subsection to show that the Betti numbers are as follows: $b_0(M_{ij})=b_3(M_{ij})=b_5(M_{ij})=b_6(M_{ij})=0$, $b_4(M_{ij})=j$ and $b_4(M_{(n-1)n})=n$. In each step, we are losing a second cohomology element, so $b_2(M_{(n-1)n})=\frac{n(n-1)}{2}$. The manifolds $M_{1j}$ are simply-connected by the Van-Kampen's theorem. As the next proposition shows, starting with $M_{23}$, the manifolds $M_{ij}$ are not simply-connected since the operation used is a self-plumbing in each step and this operation increases the first Betti number by one. Consequently, $b_1(M_{(n-1)n})=\frac{(n-1)(n-2)}{2}$. 

{\proposition Let $W'$ be a $6$-manifold obtained by plumbing two embedded spheres $a$ and $b$ in a $6$-manifold $W$ as above. Then $b_1(W')=b_1(W)+1$.}
\begin{proof} We make the attachment over $\amalg_{2}S^2$, the disjoint union of $2$ spheres. If we remove the plumbed sphere from $W'$, what we have is $W$ with the two spheres removed. This manifold is connected and the last part of the Mayer-Vietoris sequence of this operation is given below.

\begin{equation} \label{MV4}
\cdots\rightarrow H_1(\amalg_{2}S^2\times S^3;\Z)\rightarrow H_1(W';\Z)\oplus H_1(S^2\times B^4;\Z) \rightarrow H_1(W-\amalg_2 S^2\times B^4;\Z)
\end{equation}

\begin{equation*}
\rightarrow \widetilde{H}_0(\amalg_{2}S^2\times S^3;\Z) \rightarrow \widetilde{H}_0(W';\Z)\oplus \widetilde{H}_0(S^2\times B^4;\Z) \rightarrow \widetilde{H}_0(W-\amalg_2 S^2\times B^4;\Z)\rightarrow 0
\end{equation*}

We get the following short sequence (\ref{MV5}) by placing the known values to the sequence (\ref{MV4}) above.

\begin{equation} \label{MV5}
0\rightarrow H_1(W';\Z)\oplus H_1(S^2\times B^4;\Z) \rightarrow H_1(W-\amalg_2 S^2\times B^4;\Z) \rightarrow \widetilde{H}_0(\amalg_{2}S^2\times S^3;\Z)\rightarrow 0
\end{equation}

Therefore, $b_1(W')=b_1(W)+1$.
\end{proof}

The generator coming from the increase of the first Betti number is also one of the generators of the fundamental group of the plumbed manifold. If we remove the interior of the normal neighborhood of the plumbed sphere from the plumbed manifold, by the Van Kampen's theorem, the fundamental group remains the same. This manifold is obtained from the unplumbed manifold by adding a $S^2\times S^3\times I$. This last one is simply-connected so no relation is introduced in this operation. Hence the fundamental group gains a free part as a consequence of the plumbing operation. 

The calculation of the homology of the boundary is done using the reduced relative homology exact sequence of the pair $(M_{ij},\partial M_{ij})$. As an instance, sequence (\ref{RMV2}) which is given below is the exact sequence that is associated with one of the plumbings producing $M_{13}$.

\begin{equation*}
0\rightarrow H_5(M_{13};\Z)\rightarrow H_5(M_{13},\partial M_{13};\Z)\rightarrow H_4(\partial M_{13};\Z)
\end{equation*}

\begin{equation*}
\rightarrow H_4(M_{13};\Z)\rightarrow H_4(M_{13},\partial M_{13};\Z)\rightarrow H_3(\partial M_{13};\Z)
\end{equation*}

\begin{equation} \label{RMV2}
\rightarrow H_3(M_{13};\Z)\rightarrow H_3(M_{13},\partial M_{13};\Z)\rightarrow H_2(\partial M_{13};\Z)
\end{equation}

\begin{equation*}
\rightarrow H_2(M_{13};\Z)\rightarrow H_2(M_{13},\partial M_{13};\Z)\rightarrow H_1(\partial M_{13};\Z)
\end{equation*}

\begin{equation*}
\rightarrow H_1(M_{13};\Z)\rightarrow H_1(M_{13},\partial M_{13};\Z)\rightarrow \widetilde{H}_0(\partial M_{13};\Z)=0
\end{equation*}

By using Poincare-Lefschetz duality for the compact manifolds with boundary, and placing the known values to sequence (\ref{RMV2}), we find the Betti numbers of $\partial M_{13}$ as $b_0=b_1=b_4=b_5=1$ and $b_2=b_3=b_2(X_1)+b_2(X_2)+b_2(X_3)-6$.

In exact sequence (\ref{RMV2}), $H_1(M_{13},\partial M_{13};\Z)$ is trivial, which means that the first homology element is killed inside the $6$-manifold. An induction argument, starting with $M_{13}$, shows that all of the first and second homology elements of the boundary are inherited from the manifold itself. Hence, considering the classification of simply-connected, spin $5$-manifolds (\cite{Sm-62}), all second and third homology elements in the boundary are coming from a bunch of $S^2\times S^3$'s that are connected summed to each other.

For $\partial M_{13}$, the boundary of the first nonsimply-connected manifold we produced, the only remaining step is to find a manifold which is connected summed to these $S^2\times S^3$'s. Let us call this manifold $N$. As before, for the self-plumbing, a new first homology generator is introduced and it contributes a free generator to the fundamental group of the boundary. Thus, the manifold $N$ has an infinite cyclic fundamental group. By a theorem of Browder (\cite{Br-66}), we can split this manifold as $Y\bigcup_{X\times S^0}(X\times I)$, where $X$ and $Y$ are simply-connected manifolds in dimension $4$ and $5$, respectively. Note that, $Y$ is a cobordism of $X$ with itself. $\widetilde{H}_i(Y)=0$ for $i\neq 4$. If $Y$ had a nonzero third homology element, $X\times I$ would contribute a fourth homology element to $N$ which is bounded by the third homology element of $Y$. This element could only be a multiple of the fundamental class which is already in the fourth homology of $Y$ by the inclusion of the boundary. Consequently, $Y$ has trivial third homology. If $Y$ had a nonzero second homology element, this element would be an element in the second homology of $X$ and hence an element of $X\times I$. After gluing, this would be a nontrivial element in the second homology group of $N$. However, $N$ has no nontrivial second homology elements. Hence by the Hurewicz theorem, $Y$ is homotopic to $X$, i. e. $Y$ is an $h$-cobordism. All obstructions to the $h$-Cobordism theorem lies in the degree two and three homology in dimension $5$. Therefore, $Y$ is homotopy equivalent to $X\times I$ and the manifold $N$ that we are after is homotopy equivalent to $S^1\times S^4$. By a result of Shaneson for the smooth, closed, orientable $5$-manifolds with infinite cyclic fundamental group (\cite{Sh-68} Theorem in p.297), this manifold is diffeomorphic to $S^1\times S^4$. In each plumbing, this new $5$-manifold must have its first and a fourth homology groups isomorphic to $\Z$. During the plumbing, the boundary changes away from those $S^1\times S^4$'s introduced before this step. Therefore, in each step an $S^2\times S^3$ is replaced by an $S^1\times S^4$.

The $S^4$ in the boundary, formed during the plumbing, bounds a fifth homology element in the $6$-manifold, hence it dies in the $6$-manifold we constructed. Closing the boundary with the boundary sum of $B^3\times S^3$'s and $B^2\times S^4$'s gives a closed manifold $M$, whose Betti numbers are given by $b_0(M)=b_6(M)=1$, $b_2(M)=b_4(M)=n$, $b_1(M)=b_5(M)=b_3(M)=0$. 

Since $\alpha_{ii}$ is characteristic, $M$ must be spin as shown in the following proposition.

{\proposition \label{spinb=n} $M$ is spin.}
\begin{proof} The proof is similar to the proof of Proposition~\ref{spinb=2}, but in this case there are $n$ homology classes on which $w_2(M)$ must be evaluated. All of these cohomology classes $x_i$ are induced by the inclusion of $X_i$ into $M$. Pulling back everything into respective $X_i$, it turns out that $w_2(M)x_i\equiv \Sigma_j w_2(X_j)x_i+\Sigma_jw_2(\eta_j)x_i$. Remember that $\eta_i$ is the $2$-disk bundle over the $4$-manifold $X_i$. In this sum, $w_2(X_j)x_i$ and $w_2(\eta_j)x_i$ $(i\neq j)$ can be omitted due to the fact that $X_i$ overlaps $X_j$ only on $\alpha_{ij}$ which is already included in $x_i$. As a result, $w_2(X_i)x_i$ and $w_2(\eta_2)x_i$ are already included in $w_2(X_1)x_i$ and $w_2(\eta_1)x_i$, respectively. Thus, $w_2(M)x_i\equiv w_2(X_i)x_i+w_2(\eta_i)x_i\equiv w_2(X_i)x_i+Q_i(\alpha_{ii},\alpha_{ii}) \equiv 0$.
\end{proof}

The intersection form of $M$ is given in the following proposition.

{\proposition Let $M$ be the $6$-manifold constructed above with $b_2=n$ and let $\mu$ be its trilinear form. Then $\mu(x_i,x_j,x_k)=Q_i(\alpha_{ij},\alpha_{ik})$.}
\begin{proof} The proof is essentially the same as in Subsection~\ref{b=2}. The intersection $\mu(x_i,x_j,x_k)$ can be written as $PD(x_i) \pitchfork PD(x_j) \pitchfork PD(x_k)$ which is equal to $X_i\pitchfork X_j\pitchfork X_k$. The embedded copy of $X_i$ intersects the embedded copy of $X_j$ on $\alpha_{ij}$. It also intersects $X_k$ on $\alpha_{ik}$. Therefore these three intersect each other on the projection of $\alpha_{ij}$ on $\alpha_{ik}$ in $X_i$ which is nothing but $Q_i(\alpha_{ij},\alpha_{ik})$.
\end{proof}

This last proposition reflects the fact that we must choose $X_i$ and the cohomology classes such that $Q_i(\alpha_{ij},\alpha_{ik})=Q_j(\alpha_{ji},\alpha_{jk})$ for all $i$, $j$, $k$. 

The first Pontrjagin class acts on the second cohomology as follows.

{\proposition $p_1(M)x_i=3\sigma(X_i)+Q_i(\alpha_{ii},\alpha_{ii})$.}
\begin{proof} The proof is similar to the proof of Proposition~\ref{pontrjagin2} in Subsection~\ref{b=2}. We may write $p_1(M)x_i=p_1(M)PD(x_i)=p_1(M)X_i$. The tangent bundle of $M$ is restricted to the embedded copy of $X_i$ as the direct sum of the tangent bundle $TX$ of $X$ and the normal neighborhood $\nu X_i|_M$ of $X_i$ in $M$. The bundle $\nu X_i|_M$ is the $B^2$-bundle over $X_i$ with Euler class $\alpha_{ii}$. Therefore, as we have seen in Subsection~\ref{b=1} (the paragragh before Example~\ref{hcp3}), $p_1(M)X_i=p_1(X_i)[X_i]+\alpha_{ii}\cup\alpha_{ii}[X_i]=3\sigma(X_i)+Q_i(\alpha_{ii},\alpha_{ii})$.
\end{proof}

We have proved the main theorem of this paper:

{\main \label{maintheorem} Let $V$ be a smooth, closed, simply-connected, spin $6$-manifold with torsion-free homology and $b_3(V)=0$. Suppose that $H^2(V;\Z)$ is isomorphic to the direct sum of $n$ copies of $\Z$, each of which is generated by $x_i$. Also suppose that $p_1(V)x_i=24k_i+4\mu(x_i,x_i,x_i)$, where $p_1$ is the first Pontrjagin class of $V$ and $\mu$ is the symmetric trilinear form of $V$.
Then we can find a collection $\{X_i\}$ of smooth, closed, simply-connected $4$-manifolds with odd intersection forms $Q_i$ and second cohomology classes $\alpha_{ij}\in H^2(X_i;\Z)$ ($1\leq i,j\leq n$) satisfying 

\begin{enumerate} \item the signature of $X_i$ is $8k_i+\mu(x_i,x_i,x_i)$,
\item $\alpha_{ii}$ are primitive, characteristic and $Q_i(\alpha_{ii},\alpha_{ii})=\mu(x_i,x_i,x_i)$,
\item if $i\neq j$ then $\alpha_{ij}$ has a smooth sphere representative in $X_i$,
\item $Q_i(\alpha_{ij},\alpha_{ik})=\mu(x_i,x_j,x_k)$,
\end{enumerate}

so that the manifold $M$ that is constructed by closing the boundaries of the plumbed $2$-disk bundles over these $4$-manifolds $X_i$ with Euler class $\alpha_{ii}$ is diffeomorphic to $V$. The plumbing of the respective bundles over $X_i$ is done on the sphere representatives of $\alpha_{ij}$ in $X_i$ and $\alpha_{ji}$ in $X_j$.}
 
\section{The Building Blocks} \label{4man}

In this section, we prove the fact that $4$-manifolds that are used as the building blocks of the constructions in Section~\ref{constructions} exist. We state two existence results. The first lemma is a special case of the second, however we include it here because it makes reading the proof of the latter lemma easier.

{\lemma \label{8k+mExistence} For all $(k,m)\in\Z\oplus\Z$, there exists a closed, simply-connected, smooth $4$-manifold $X$ with an odd intersection form $Q$ and a primitive characteristic cohomology class $\alpha\in H^2(X;\Z)$ such that $Q(\alpha,\alpha)=m$ and the signature, $\sigma(X)$, of $X$ is $8k+m$.}
\begin{proof}
The strategy of the proof is to find the manifolds for $(k,0)$ and $(0,m)$ and then give the remaining ones as the connected sums of these manifolds. Note that, the signature of a connected sum is equal to the sum of signatures of the components. The characteristic class $\alpha_{k,m}$ of the resulting manifold after connected summing is taken as the sum of the characteristic classes of each component manifold. 
\begin{enumerate} \item $(k,m)=(0,0)$: Take $X_{0,0}=\CP\#\CPB$ and $\alpha_{0,0}=h+e$, where $h$ is the generator of $H^2(\CP;\Z)$ and $e$ is the generator of $H^2(\CPB;\Z)$. Then $Q(\alpha_{0,0},\alpha_{0,0})=0=m$ and $\sigma(X_{0,0})=8k+m=0$.
\item \begin{enumerate} \item $m=0$, $k>0$: Assume that $X_{k,m}=\#_{16k+1}\CP\#_{8k+1}\CPB$ and $\alpha_{k,0}=\sum_{i=1}^{16k+1}h_i + \sum_{i=1}^{k}3e_i +\sum_{k+1}^{8k+1}e_i$, where $h_i$ is the generator of the second cohomology group of the $i$th copy of \CP and $e_i$ is the generator of the second cohomology group of the $i$th copy of \CPB. Then the signature of $X_{k,0}$ is $8k=8k+m$ and $Q(\alpha_{k,m},\alpha_{k,m})=0=m$. 
\item $m=0$, $k<0$: Let $X_{k,m}=\#_{-8k+1}\CP\#_{-16k+1}\CPB$ and $\alpha_{k,m}=\sum_{i=1}^{-16k+1}e_i + \sum_{i=1}^{-k}3h_i +\sum_{-k+1}^{-8k+1}h_i$, where $h_i$ is the generator of the second cohomology group of the $i$th copy of \CP and $e_i$ is the generator of the second cohomology group of the $i$th copy of \CPB. Then the signature of $X_{k,m}$ is $8k=8k+m$ and $Q(\alpha_{k,m},\alpha_{k,m})=0=m$. 
\end{enumerate}
\item \begin{enumerate} \item $k=0$, $m>0$: Take $X_{k,m}=\#_{m}\CP$ and $\alpha_{k,m}=\sum_{i=1}^{m}h_i$, where $h_i$ is the generator of the second cohomology group of the $i$th copy of \CP. The signature of $X_{k,m}$ is $m$ and $Q(\alpha_{k,m},\alpha_{k,m})=m$.
\item $k=0$, $m<0$:  Take $X_{k,m}=\#_{-m}\CPB$ and $\alpha_{k,m}=\sum_{i=1}^{-m}e_i$, where $e_i$ is the generator of the second cohomology group of the $i$th copy of \CPB. The signature of $X_{k,m}$ is $m$ and $Q(\alpha_{k,m},\alpha_{k,m})=m$.
\end{enumerate}
\item For arbitrary $(k,m)$, let $X_{k,m}$ be chosen as the connected sum of $X_{k,0}$ and $X_{0,m}$ and $\alpha_{k,m}=\alpha_{k,0}+\alpha_{0,m}$. 
\end{enumerate} 
\end{proof}

Our second existence result is given by the following lemma.

{\lemma \label{b2=nExistence} Given $(k,n)\in\Z\oplus\Z_+$ and a symmetric matrix $A=\{A_{ij}\}\in GL(n;\Z)$ it is always possible to find a closed, simply-connected, smooth $4$-manifold $X$ with an odd intersection form $Q$ and $n$ distinct primitive second cohomology classes of $X$ satisfying the following conditions:
\begin{enumerate} \item one of the cohomology classes, $\alpha=\alpha_{1}$, is characteristic, and $Q(\alpha,\alpha)=A_{11}$,
\item the $(n-1)$ remaining cohomology classes $\{\alpha_i\}_{i=2}^{n-1}$ are represented by embedded spheres, 
\item the intersection numbers of the $n$ cohomology classes are given by $Q(\alpha_i,\alpha_j)=A_{ij}$,
\item the signature, $\sigma(X)$, of $X$ is equal to $8k+Q(\alpha,\alpha)$.
\end{enumerate}}

\begin{proof} When $n=1$, the proof is given in Lemma~\ref{8k+mExistence}. Let $n>1$. We give the homeomorphism type of $X$ by writing it as a connected sum of \CP's and \CPB's. There is a distinguished class $\alpha$ which is characteristic. For each $\alpha_i$, we find a manifold $X_i$ corresponding to $\alpha_i$. Then we glue all $X_i$ together according to the intersection of $\alpha_i$ with the other classes. We start gathering the pieces of $X$ by considering the difference between $Q(\alpha_i,\alpha_i)$ and $Q(\alpha_i, \alpha)$. This number is always even (since $\alpha$ is characteristic). Let $Q(\alpha_i,\alpha)-Q(\alpha_i, \alpha_i)$ be equal to $2k_i$. The manifold $X_i$ and the class $\alpha$ are formed by the following algorithm:
\begin{enumerate}\item This step is to adjust the intersection with the characteristic class $\alpha$. If $k_i$ is positive, then $X_i=\CP$, $\alpha_i=h_i$ and $\alpha=(2k_i+1)h_i$, where $h_i$ is the generator of the second cohomology group of \CP. If $k_i$ is negative, then $X_i=\CPB$, $\alpha=(2k_i+1)e$ and $\alpha_i=e_i$. Here $e_i$ is the generator of the second cohomology of \CPB.
\item Now we adjust the intersection number with other cohomology classes. If $Q(\alpha_i,\alpha_j)$ is negative, $X_i$ and $\alpha_i$ stay as they are. If $Q(\alpha_i,\alpha_j)$ is positive then add $Q(\alpha_i,\alpha_j)$ copies of \CP's to $X_i$ by connected summing, and add the generators of second cohomologies of the new \CP's to $\alpha_i$. The same number of \CP's are also added to $X_j$ but $\alpha_j$ is not changed. 
\item  To reach the self-intersection  $Q(\alpha_i,\alpha_i)$ of $\alpha_i$, connect sum $X_i$ with a number of \CP's or \CPB's, and add the generators of the new cohomology to $\alpha_i$. The number of the new manifolds that are added is determined by the self-intersection of $\alpha_i$ at the end of the last step. 
\end{enumerate}
Now we have a collection of manifolds $\{X_i\}$. The next step is to glue all these manifolds by identifying $X_i$ and $X_j$ on the \CP's added to $X_i$ and $X_j$ in the second step of the algorithm above. In this manifold, the characteristic class $\alpha$ is given by $\sum(2k_i+1)h_i+\sum(2k_i+1)e_i$. As a result, we have a manifold with $n$ primitive cohomology classes.The intersection numbers of these classes were adjusted except $Q(\alpha,\alpha)$. To obtain the desired manifold, we need to adjust the self-intersection of $\alpha$ and the signature. We manage the self-intersection of $\alpha$ by adding \CP's (or \CPB's) if necessary. $\alpha$ is changed by the addition of a bunch of $h$'s (or $e$'s). Finally, we reach the signature by connected summing with the manifolds used in the proof of Lemma~\ref{8k+mExistence}. The contribution of these manifolds to $\alpha$ does not change the self-intersection of $\alpha$.

The projection of $\alpha_i$ in each $\CP$ ($\CPB$) can be represented by a sphere. This is because it is nothing but a line $h$ in $\CP$ (exceptional sphere $e$ in $\CPB$). Connecting these spheres linearly by thin tubes, we can represent all the $\alpha_i$'s by embedded spheres. 
\end{proof}


\end{document}